\documentclass[3p]{elsart}
\usepackage{xcolor}
\usepackage{times}
\usepackage{helvet}
\usepackage{courier}
\usepackage[T1]{fontenc}
\usepackage{amsmath}
\usepackage{amsfonts}
\usepackage{marvosym}
\usepackage{textcomp}
\usepackage{psfrag}
\usepackage{calc}
\usepackage{graphicx}
\usepackage[normalem]{ulem}
\usepackage{algorithm}
\usepackage{latexsym}
\usepackage{multirow}

\usepackage{amsfonts}
\usepackage{marvosym}
\usepackage{psfrag}
\usepackage{algorithm}
\usepackage{pstricks,pst-node,pst-tree}
\usepackage{centernot}

\usepackage{array}
\usepackage{rotating}
\newcommand{\rot}[1]{ \begin{minipage}{6.2mm}\begin{rotate}{55}~\hspace{-2.9mm} #1 \end{rotate} \end{minipage}}

\newcommand{\ds}{\displaystyle}
  \newtheorem{dfn}{Definition}
  
  \newtheorem{example}{Example}
  \newtheorem{crt}{Criterion}
 \newtheorem{remark}{Remark}

\newcommand{\R}{\mathbb{R}}

\begin{document}
\theoremstyle{plain}

\begin{frontmatter}
\title{Diagnosability and detectability of multiple faults in nonlinear models}


\author{N. Verdi\`ere$^{1}$}
\author{S. Orange$^{1}$}

\address{$^{1}$ Normandie Univ, UNIHAVRE, LMAH, FR-CNRS-3335, ISCN, 76600 Le Havre, France\\
(Corresponding author: nathalie.verdiere@univ-lehavre.fr)}  
          

\begin{abstract}
This paper presents a novel method for assessing multiple fault diagnosability and detectability of nonlinear parametrized dynamical models. This method is based on computer algebra algorithms which return precomputed values of algebraic expressions characterizing the presence of some constant multiple fault(s). Estimations of these expressions, obtained from inputs and outputs measurements, permit then the detection and the isolation of multiple faults acting on the system. This method applied on a coupled water-tank model attests the relevance of the suggested approach.

\end{abstract}

\end{frontmatter}

\section{Introduction}

The problem of fault-diagnosis has received an increasing attention during the recent years in order to increase security of systems, to monitor their performance or to endow them with self diagnostic capabilities. To answer such technological requirements, this problem needs to be taken into account in the system design stage from an a priori diagnosability study on a model. In studying anticipated fault situations from different symptoms of the system, faults or multiple  faults can be known as discriminable according to the available sensors in a system. Some procedures for detecting and isolating them may, then, be put in place in the design stage. By this way, diagnosability can permit to anticipate component failures. 

In this paper, we assume to be in the model-based framework and, more precisely, that available signals $u$ and $y$ and a nonlinear parametri\-zed dynamical model permit to reach the output trajectories of the system. The problem consists in this framework to evaluate diagnostic performance given a model only. By (multiple) fault, we mean any change(s) of  parameter value(s) implying unwanted changes in the behavior of one or more component(s) of the system. The \textit{fault diagnosis} study consists in two subtasks \cite{Gertler}. The first one concerns the \textit{fault detection} (FD) of the malfunction, the second one the \textit{fault isolation} (FI) of the faulty component (that is the determination of its location). The fault diagnosis is done from the comparison between predictions of the model and behaviors of the system. Several methods are proposed in the literature as nonlinear observers \cite{Seliger91} or methods based on testable subsets of equations \cite{Armenglo2009}. The issue of subsets generation has been studied by many authors. They can be based on Minimal Structurally Overdetermined sets (MSO) \cite{Gelso2008,Krysander2008}, on possible conflicts \cite{Pulido2004} or on Analytical Redundancy Relations (ARRs) \cite{Trave2006}. The latter are relations linking inputs, outputs, their derivatives, the parameters of the model and the faults (See~\cite{Cruz-Victoria08,Daigle2012,Marcel,Nath2015,Zhang98} for single faults and \cite{Multiplefaults1,Multiplefaults3,Multiplefaults2} for multiple faults). Some of these ARRs are obtained using computer algebra tools such as the Rosenfeld-Groebner algorithm which permits to eliminate the unknown variables of the model. 
With respect to a specific elimination order, this algorithm returns particular differential polynomials classically called \textit{input-output polynomials} (See~\cite{Boulier:1997,dj,art1}). Some recent works have already used these particular polynomials in diagnosis assuming that the model is identifiable with respect to the faults. 
Indeed, identifiability insures that the fault values can be uniquely inferred from input-output measurements. In the case of single faults, authors in~\cite{Nath2015}  prove that if the model is identifiable with respect to the faults then all the faults are discriminable; in other words, the model is diagnosable. Furthermore, they prove that the residuals associated to each ARR permit to detect each identifiable fault in adopting a discriminable behavior. In \cite{Jauberthie13}, assuming that the faults act only additively on parameters, detectability is obtained directly from the ARRs. In this last paper, interval analysis is used to estimate the simple faults. \\

We propose a new approach to exploit such ARRs to discriminate, detect and isolate (multiple) fault(s) in models not necessary identifiable. Starting from the model, the three first step of our method described hereafter can be completely automatized; our contributions consist in the steps (2) and (3).
\begin{enumerate}
 \item The first one is the computation of ARRs by applying the Rosenfeld-Groebner algorithm to the model. 
 \item The second step consists in using Groebner basis computations in order to obtain an algebraic application called \textit{algebraic signature}. Each of its components depends only on the parameters and on the coefficients of the ARRs which can be numerically estimated from the inputs and output measurements of the system. By construction, each component of the algebraic signature vanishes when at least one specific (multiple) fault occurs.
 \item For each possible (multiple) fault, the third step consists in using semialgebraic set tools to certify that some components of the signature vanish or never vanish.  These expected values can be summarized in a precomputed table. 
 \end{enumerate}

This table constitutes the input of the numerical treatment which, from the system measurements, returns estimations of the algebraic signatures. Their comparison with their nominal values permits to detect and isolate (multiple) fault(s). \\

The advantage of the present method is, first, not to require any strong assumption on the model for determining some possible acting multiple fault(s) as i) its identifiability, ii) the value of some model parameters  in some particular cases, iii) that the multiple faults act only additively on parameters. Then, by using semialgebraic sets tools, constraints on parameters and multiple faults, such as inequalities satisfied by parameters or constraints deduced from initial conditions, can be taken into account through automatic procedures. These constraints can play a fundamental role in the FDI analysis of a model as it will be seen in this paper. Finally, from an a priori study on the model, a numerical method based only on estimation of algebraic expressions, and consequently fast, is proposed to do FDI.\\

The paper is organized as follows. In Section \ref{SectionModel}, we present the framework of our method and we define the notion of algebraic diagnosability of multiple faults. Section \ref{SecConstructionAsign} is devoted to our method consisting in studying the diagnosability of a model, that is the way to compute an algebraic signature and to tabulate its expected values in function of the multiple faults. In Section \ref{SecAppli}, our method is applied to an example of two coupled water-tanks. Section \ref{conclusion} concludes the paper. In this paper, the symbolic computations had been realized with Maple 18 and the numerical part with Scilab.

\section{Dynamical models and Diagnosability}\label{SectionModel}
\subsection{Studied parametrized models}\label{sectiondef}
We consider nonlinear parametrized models controlled or uncontrolled of the following form:
\begin{equation}\label{saresoudre}
\Gamma_f
\left\{\begin{array}{l}
\dot{x}(t,p,f)=g(x(t,p,f ),u(t),p,f ),\\
y(t,p,f)=h(x(t,p,f ),u(t),p,f ),\\
t_0\leq t\leq T
\end{array}\right.
\end{equation}
where: 
\begin{itemize}
\item the vector of real parameters $p =(p_1,\ldots,p_m)$ belongs to $\mathcal{P} \subseteq \R^{m}$ where  $\mathcal{P}$ is an \textit{a priori }known set of admissible parameters,
\item $f =(f_1,\ldots,f_e)$ is a constant fault vector which belongs to a subset $\mathcal{F}$ of $\R^e$. It is equal to 0 when there is no fault. The set $\mathcal{F}$ describes the set of admissible values of the fault vectors $f $,
\item $x(t,p,f )\in \R^n$ denotes the state variables and $y(t,p,f )\in \R^s$ the outputs,
\item  $g$ and $h$ are real vectors of rational functions in $x$, $p$ and $f$ \footnote{This assumption is not restrictive since lots of models can be reduced to a rational model by variable change (see.~\cite{SAB}).}.
\item $u(t)\in \R^r$ is the control vector equal to 0 in the case of uncontrolled models.
\end{itemize}

\begin{remark}
In most practical cases, the faults $f_i$ belong to connected sets of $\mathbb{R}$, and $\mathcal{F}$ is the Cartesian product of these sets. The present work takes place in a more general framework by introducing semialgebraic sets defined hereafter.
\end{remark}
From now on, we suppose that constraints on $p \in \mathcal{P} $ and $f \in \mathcal{F}$, and eventual constraints linking faults and parameters components, can be formulated by the mean of algebraic equations and/or inequalities. This leads naturally to consider semialgebraic sets for which computer algebra tools are developed (See~\cite{qepcadb,RAGLib,DISCOVERED}, for example):

\begin{dfn}\label{defsemialgebraicset} {(See~\cite{Basu})} A set of real solutions of a finite set of polynomial equations and/or polynomial inequalities is called a semialgebraic set.
\end{dfn}

Let $C_{p,f}$ be the set of all algebraic equations and inequalities verified by the components of the parameter and fault vectors of the model and $\mathcal{C}_{p,f}$ be the semialgebraic set defined by $C_{p,f}$. In order to take into account initial conditions, the algebraic relations induced by these conditions can be added to the set $C_{p,f}$. 


Let $\mathcal{N}$ be a subset of $\{1,\, \ldots, \, e\}$ and $f_\mathcal{N}$ the \textit{multiple  fault vector} whose components $f_i$ are not equal to 0 if $i \in \mathcal{N}$ and equal to 0 otherwise. Naturally, $f_\mathcal{N}$ belongs to $ \mathcal{F}_\mathcal{N} = \{ f \in \mathcal{F} | f_i \neq 0\,\text{ if }\, i \in \mathcal{N}\,\text{ and }\, f_i = 0\,\text{ if }\, i \notin \mathcal{N}\, \} $ and $\mathcal{F}_\mathcal{N}$ is a semialgebraic set by construction.
When only one component of $f$ is not null, the fault vector $f$ is called a \emph{simple fault}.

Since diagnosability and detectability of Model~(\ref{saresoudre}) may depend on $p =(p_1,\ldots,p_m)$, we consider, afterwards, the set $R=\R[p_1,\ldots,p_m]$ of polynomials in the indeterminates $p_i$ with real coefficients.\\

\subsection{Algebraic signature and diagnosability}\label{sectionAS}

To characterize multiple faults, the following definition introduces the notion of algebraic signature. It is based on $l$ algebraic expressions, $ASig_i$ ($i=1,\ldots,l$), deduced from the system~$\Gamma_f$. These expressions depend on the model parameters and the fault components.

\begin{dfn}
Let $ASig=(ASig_1,\, \ldots ,\, ASig_l )$ be a vector of polynomial functions admitting $f_1,\; \cdots,\; f_e$ as indeterminates with coefficients in $ {R}$. An algebraic signature is a function $ASig$ defined by: \\
\begin{center}
$\begin{array}{cccc} \label{functionasign}
ASig : & \mathcal{F} &\longrightarrow &R^{l} \\
 & f &  \mapsto & (ASig_1(f),\, \ldots,\, ASig_l(f) ). 
\end{array}$\end{center}
\end{dfn}

The comparison of the image of two multiple faults under the function $ASig$ gives a way to discriminate them.
If the model is controlled, we propose to define the strongly and weakly algebraic diagnosability, the first one being true for all inputs and the second one for at least one input. 
\begin{dfn}\label{strongly_fd}
Let $\mathcal{N}$ and $\mathcal{N}'$ be two distinct subsets of $\{1,\, \ldots, \, e\}$. The multiple faults of $\mathcal{F}_\mathcal{N}$ and of $ \mathcal{F}_{\mathcal{N}'}$ are said input-strongly algebraically discriminable (resp. input-weakly algebraically discriminable) if, for all input $u$ (resp. one input),
\begin{equation}\label{isfd}
\begin{array}{l}
ASig(\mathcal{F}_{\mathcal{N}})\cap ASig(\mathcal{F}_{\mathcal{N}'}) = \emptyset .
\end{array}
\end{equation}
 This equality is in particular satisfied when there exists an index $i$ such that $ASig_i(\mathcal{F}_\mathcal{N}) \cap ASig_i(\mathcal{F}_{\mathcal{N}'})=\emptyset\;. $

If, for any distinct subsets $\mathcal{N}$ and ${\mathcal{N}'}$ of $\{1,\, \ldots, \, e\}$, the multiple faults of $\mathcal{F}_\mathcal{N}$ and $ \mathcal{F}_{\mathcal{N}'}$ are input-strongly algebraically discriminable (resp. input-weakly algebraically discriminable), the model is said input-strongly algebraically diagnosable (resp. input-weakly algebraically diagnosable).
\end{dfn}
In the case of uncontrolled model, the definition of algebraic diagnosability can be proposed too in omitting the notion of input in the previous definitions.\\




Naturally, the notion of detectability of a set of multiple faults can be defined from the algebraic signature in comparing its value to the one obtained when the set $\mathcal{N}$ is empty, that is when no fault occurs in the system. This definition is given below. 
\begin{dfn}
A set of multiple faults vectors $\mathcal{F}_\mathcal{N}$ is algebraically detectable if $$ASig(\mathcal{F}_\mathcal{N}) \cap ASig(\mathcal{F}_\mathcal{\emptyset})= \emptyset ,\, $$
$ASig(\mathcal{F}_\mathcal{\emptyset})$ being the algebraic signature evaluated when no fault occurs in the system.
\end{dfn} 

In Section~\ref{SecConstructionAsign}, we propose a procedure to obtain an algebraic signature not depending explicitly on the (multiple) faults and assessable from the known quantities of the system. Some criterions are also proposed to discriminate as far as possible all the (multiple) faults. 

\section{An algebraic diagnosability method} \label{SecConstructionAsign}
For practical applications, we propose the construction of an algebraic signature in three steps using symbolic computations. The first one requires the implementation of the Rosenfeld-Groebner algorithm in order to obtain algebraic relations linking parameters, faults and real values deduced from the outputs of the system. From these algebraic relations, the second step consists in using the Groebner basis algorithm to obtain an algebraic signature and the third one uses semialgebraic set tools to discriminate multiple faults.

\subsection{First step: construction of the exhaustive summary from the model}
In \cite{Nath2015}, the authors give a way to obtain relations linking inputs, outputs, parameters and faults. The latter are obtained from the Rosenfeld-Groebner algorithm implemented in some computer algebra systems. This elimination algorithm used with an appropriate elimination order permits to eliminate unknown variables from System~(\ref{saresoudre}). These input-output representations may act as analytical redundancy relations (ARRs) and have the following forms\\
$$\ds w_i(y,u,p ,f)=m _{0,i}(y,u,p )+\sum_{k=1}^{n_i}\gamma_k^i(p,f ) m_{k,i}(y,u)\, , \\ \,i=1,\ldots,s$$
where $(\gamma_k^i)_{1\leq k\leq n_i}$ are rational fractions in $p$ and $f$, $\gamma_k^v\ne \gamma_k^w$ for $v\ne w$, $(m_{k,i}(y,u))_{1\leq k\leq n_i}$ are differential polynomials with respect to $y$ and $u$ and  $m _{0,i}\ne 0$.\\ The first part of this polynomial is supposed not to be identically equal to zero and does not contain components of $f$. It corresponds to the residual \textit{computation form} whereas the second form is known as the residual \textit{internal form}. According to \cite{dj}, there are as many polynomials of this form as outputs.  

The sequence $(\gamma_k^i(p,f ))_{k=1,\ldots,n_i}$ ($i=1,\ldots,s$) is called the exhaustive summary of System~(\ref{saresoudre}) (See~\cite{NOLCOS}). We now consider the function $\phi$ constructed from the exhaustive summary defined by:
\begin{center}
$\begin{array}{cccc} \label{deffunctionphi}
\phi : &  \mathbb{R}^e &\longrightarrow &R^{N } \\
 & f &  \mapsto & (\gamma _k^i(p,f ))_{1\leq i\leq s,1\leq k\leq n_i}
\end{array}$\end{center}
where $\displaystyle N=\sum_{i=1}^{s} n_i$. 

To lighten our approach, we suppose that $(\gamma_k^i)_{1\leq k\leq n_i}$ are polynomials of $R [f_1,\ldots,f_e]$ \footnote{Actually, when these expressions are rational fractions, non vanishing conditions for the denominators can be added to $C_{p,f} $ and new variables corresponding to the inverse of the denominators can be added to rewrite $(\gamma_k^i)_{1\leq k\leq n_i}$ as polynomials.} where $R=\R[p_1,\ldots,p_m]$.

\begin{example}\label{exemple1}
Consider the ARR
$$y^2+(p_1+f_1)\dot y^2+p_2( (f_2-1)^2-1/4)y\dot y + p_3 f_1 - p_1=0$$
where $p=(p_1,p_2,p_3)$ and $f=(f_1,f_2)$ have respectively their components in $]0,+\infty )$ and in $[0,2[$. \\
The function $\phi$ is then defined by \begin{center}
$\phi(f)=(p_1+f_1,p_2((f_2-1)^2-1/4), p_3 f_1 - p_1)$.                                                                               \end{center}
\end{example}


By definition, $\phi(f)$ defines an algebraic signature. The injectivity of $\phi$ is strongly connected to the notion of identifiability of the model. Recall that a model is identifiable if the model parameters are uniquely determined by the model inputs and outputs. In \cite{NOLCOS}, the authors prove that if the function $\phi$ is injective and under some technical assumptions, the model is identifiable. Consequently, full identifiability of the fault parameters implies algebraic diagnosability since any fault vector instance will give a distinct value of $\phi(f)$. However, algebraic diagnosability does not imply identifiability. Indeed, even if $\phi$ is not injective, (multiple) faults discrimination may be possible as shown  in the following example.
\begin{example}
In example \ref{exemple1}, the function $\phi$ is not injective: the values $1/2$ and $3/2$ of $f_2$ will give the same value of $\phi(f)$.\\ By setting $ASig = \phi$, the algebraic signatures of the possible multiple faults are $ASig(f_{\emptyset}) =(p_1,\frac{3}{4}p_2 )$, $ASig(f_{\{1\}})=(p_1+f_1,\frac{3}{4}p_2)$, $ASig(f_{\{2\}})=(p_1, p_2((f_2-1)^2-\frac{1}{4}))$ and $ASig(f_{\{1,2\}})=(p_1+f_1, p_2((f_2-1)^2-\frac{1}{4}))$. Constraints on parameters and faults imply that the intersection of the images of these algebraic signatures do not intersect. Consequently, the model is algebraic diagnosable since the multiple faults can be discriminated.
\end{example}
  The algebraic signature defined by the exhaustive summary is not sufficient since two distinct faults acting on its same components may not be discriminated. A natural approach to exploit the exhaustive summary consists in obtaining an explicit expression of the fault components in function of $\phi$ and the model parameters. This approach focusing on the inversion of an algebraic system fails in general. That is why we propose a method to obtain algebraic expressions not depending on the faults and characterizing their presence. Such expressions can be generated by automatic procedures based on Groebner basis computations (See~\cite{cox, Faugere}) and are used, in the next section, to define an algebraic signature.   

\subsection{Second step: construction of an algebraic signature from the exhaustive summary}

Given a multiple fault $f\in \mathcal{F}_\mathcal{N}$ ($\mathcal{N} \subset \{1,\; \cdots,\; e\}$), let $E_\mathcal{N}$ be the set of polynomials\\
$\begin{array}{l}
E_\mathcal{N} = \{\gamma _1^1(p,f ) - \phi_1,\; \ldots,\; \gamma _s^{n_s}(p,f) - \phi_N\}\, \\
~~~~~~~~~~~~~~~~~~~~~~\cup \,\{ v_i f_i -1 | i \in \mathcal{N} \}\, \cup \, \{  f_i  | i \notin \mathcal{N} \}\, \end{array}$\\
where $v_i$ are new indeterminates. 
In the definition of $E_\mathcal{N}$, the sets $\{ v_i f_i -1 | i \in \mathcal{N} \}$ and $\{  f_i  | i \notin \mathcal{N} \}$ characterize multiple  faults of $\mathcal{F}_{\mathcal{N}}$. Let us consider the polynomial ideal $I_{\mathcal{N}}$ generated by $E_\mathcal{N}$, that is the set of all linear combinations of elements of $E_{\mathcal{N}}$ in $R [v_1,\ldots,v_e, f_1,\ldots,f_e, \phi_1,\ldots, \phi_N]$. 

A Groebner basis of this ideal $I_{\mathcal{N}}$ is computed  with respect to an elimination order chosen to eliminate first the indeterminates $v_i$ and  $f_i$. The intersection $G_{\mathcal{N}}$ of this Groebner basis and of $ R [\phi_1,\ldots, \phi_N]$ generates the elimination ideal  $J_\mathcal{N}= I_\mathcal{N}\;  \cap \;R [\phi_1,\ldots, \phi_N]$ (See~\cite{cox}). Clearly, any polynomial of $G_\mathcal{N}$ vanishes when a multiple  fault $f\in \mathcal{F}_\mathcal{N}$ occurs.\\
For all the possible multiple faults $f_\mathcal{N}$, the sets $G_\mathcal{N}$ are computed. Polynomials of $\cup_{\mathcal{N}\subset\{1,\, \ldots,\,m\} }G_\mathcal{N}$ vanishing for all multiple faults, i.e. polynomials of $\cap_{\mathcal{N}\subset\{1,\, \ldots,\,m\}} I_\mathcal{N} $, are removed of this set. The remaining polynomials are kept to define the components of an algebraic signature.

Let us summarize our algorithm returning an algebraic signature.

\verb Algebraic_signature  \hrulefill

\begin{enumerate}
 \item For each subset $\mathcal{N}$ of $\{1,\; \cdots,\; e\}$, we consider a generic multiple fault $f_\mathcal{N}$ and we apply the following steps to this multiple  fault.
 \begin{enumerate}
  \item Computation of the Groebner basis of the ideal $I_\mathcal{N}$ generated by $E_\mathcal{N}$ with respect to the lexicographical order $v_{i_1} \succ \ldots  \succ  v_{i_l} \succ f_1\succ \ldots\succ f_m \succ  \phi_1 \succ \ldots \succ  \phi_N \succ p_1 \succ \, \dots \succ p_m .$ 
  \item Determination of the intersection, $G_\mathcal{N}$, of this last Groebner basis and of $R [ \phi_1,\ldots, \phi_N]$. 

 \end{enumerate}
 \item Remove to $\cup_{\mathcal{N}\subset\{1,\, \ldots,\,m\}}G_\mathcal{N}$ polynomials vanishing for any multiple fault, in other words, polynomials of the ideal $\cap_{\mathcal{N}\subset\{1,\, \ldots,\,m\}} I_\mathcal{N} $. 
 \item Order arbitrarily all the polynomials of the last obtained set in a sequence $ASig=(ASig_1,\, \ldots,\, ASig_l )$.
 \item Return $ASig$. \end{enumerate}
 \hrulefill
 
By this way, we obtain an \textit{algebraic signature of the multiple faults} used afterwards:\\
$$\begin{array}{cccc} \label{deffunctionASig}
ASig : &  \mathbb{R}^e  &\longrightarrow &(R[ \phi_1,\ldots, \phi_N])^{l} \\
 & f &  \mapsto & (ASig_1( \phi),\, \ldots,\, ASig_l( \phi) ) \, . 
\end{array}$$

\begin{example}\label{exemple3}
Let us consider the exhaustive summary  $  \phi (f_1, f_2)= (p_1+f_1, p_2((f_2-1)^2-\frac{1}{4}),p_3 f_1 - p_1)$ of Exemple~\ref{exemple1}. 
The algorithm \verb Algebraic_signature   returns the signature $ASig$ defined by $ASig(f_1,f_2) = (\phi_1-p_1, 4\phi_2-3p_2 , p_1+\phi_3 )$ whose components vanish for at least one mutiple fault.
\end{example}

By construction, the signature $ASig(f)$ does not depend explicitly on $f$. However, the presence of multiple  fault(s) is reflected in the numerical values of $\phi$ and, consequently, of $ASig$. 
From the comparison between an estimation of $ASig(f)$ and the expected null components of the lists $ASig(f_\mathcal{N})$, some possible multiple  faults can be discarded. 
Nevertheless, such a comparison may not be sufficient to discriminate some multiple  fault signatures. 
Indeed, polynomials of $G_\mathcal{N}$ appearing in $ASig(f_\mathcal{N})$ are insured to vanish when the fault $f_\mathcal{N}$ occurs but the other components of the signature $ASig(f_\mathcal{N})$ may also vanish for some particular values of the parameters and faults. That is why it is necessary to introduce supplementary criterions to improve the multiple faults discrimination.

\subsection{Third step : Criterions to differentiate multiple fault signatures}\label{SectionCriterions}

In order to elaborate additional criterions, the\textit{ semialgebraic approach} (See~\cite{Basu}), focusing on real solutions of polynomial equations and inequalities, is adapted. This approach permits to take into account the set of constraints on parameters and on faults, $C_{p,f}$, of System~(\ref{saresoudre}); this set can play an important role for the discrimination of multiple fault signatures as explained in Exemple~\ref{exemple3suite}.

The three following results lies on the emptyness of semialgebraic sets which can be tested by using computer algebra tools (See~\cite{RAGLib,DISCOVERED}). The first criterion (resp. the second) consists in determining whether the $k$-th component of $ASig(f_\mathcal{N})$ vanishes for at least one real values of a multiple  fault $f \in \mathcal{F}_\mathcal{N}$ (resp. never vanishes).

 For any $\mathcal{N} \subset \{1,\; \cdots,\; m\}$, let us consider the set $S_\mathcal{N}$ of polynomial equations and inequalities defined by \\
$S_\mathcal{N} = \{\gamma _1^1(p,f ) =  \phi_1,\; \ldots,\; \gamma _s^{n_s}(p,f) = \phi_N\} \\~~~~~~~~~~~~~~\,\cup\, C_{p,f} \, 
\cup \,\{ v_i f_i =1 | i \in \mathcal{N} \}\,\cup \, \{  f_i =0  | i \notin \mathcal{N} \}.$
\begin{crt} If the semialgebraic set defined by $S_\mathcal{N}\, \cup  \{ ASig_k(f_{\mathcal{N}})=0\}$ is empty then the $k$th component of $ASig(f_\mathcal{N})$ never vanishes.\label{Criterion1}
\end{crt}

\begin{crt} If the semialgebraic set defined by $S_\mathcal{N}\, \cup \{ ASig_k(f_{\mathcal{N}}) v_k-1 =0\}$ is empty then the $k$th component of $ASig(f_\mathcal{N})$ is equal to 0.\label{Criterion2}
\end{crt}

For some particular system, a vanishing component of the signature may charaterize multiple  faults $f$ whose $i$th component is not null.

\begin{crt} Let $S$ be the semialgebraic set defined by $$ S = \{\gamma _1^1(p,f )=  \phi_1,\; \ldots,\; \gamma _s^{n_s}(p,f) = \phi_N\}\,\cup\, C_{p,f}.$$\label{Criterion3}
If the sets of real solutions $S \, \cup\, \{ASig_j(f)=0,v_i\, f_i-1 =0\}$ and $S \, \cup\, \{v_j\,ASig_j(f)-1=0,f_i=0\}$ are empty then $ ASig_j(f)=0$ is equivalent to $f_i = 0$.
\end{crt} 

In the case where Criterion~\ref{Criterion3} is satisfied for all the components $f_i$ of $f=(f_1,\, \ldots, \, f_m)$, it is clearly useless to apply criterions~\ref{Criterion1} and ~\ref{Criterion2} on all the $m!$ possible multiple  faults since it permits to determine non null components of $f$.


With the help of these three criterions, the expected values of $ASig(f)$ when a multiple  faults $f$ occurs can be tabulated. In the next example and in Section~\ref{SecAppli}, the following convention is used in these tables: for any multiple fault $f$, 
\begin{itemize}
 \item A 0 in a cell means that the corresponding component of $ASig_i(f)$ is necessarily equal to 0. This is the case when $ASig_i(f)$ belongs to $I_\mathcal{N}$ and this can also be a consequence of Criterion~\ref{Criterion2}.
 \item A cell containing $\centernot{0}$ means that Criterion~\ref{Criterion1}  insures that the component of $ASig(f)$ never vanishes when the multiple  fault occurs.
 \item An empty cell signifies that the component of the signature vanishes for some values of $(p,f)$ and does not vanish for some other values of $(p,f)$.
\end{itemize}

\begin{example}\label{exemple3suite}
Let us continue Exemple~\ref{exemple1} and consider the algebraic signature function $ASig(f_1,f_2) = (\phi_1-p_1, 4\phi_2-3p_2 , p_1+\phi_3 )$ returned by the second step of our method (See Example~\ref{exemple3}).

If the set of constraints $ C_{p,f} = \{0<p_1,0<p_2,0<p_3,0\leq f_1<2,0\leq f_2<2,0\leq f_3<2\}$ is taken into account, the two first criterions  provide some characteristics of the algebraic signature for the possible mutiple faults. They are summarized in the following table:
\begin{center}
\begin{tabular}{|c|c|c|c|}
\hline  
$f$  & {$ ASig_1(f)$}  & {$ ASig_2(f)$} & {$ ASig_3(f)$} \cr 
      \hline $f_{\{\}}$ &$0$&$0$&$0$\cr
     \hline {$f_{\{1\}}$} & $\centernot{0}$ &$0$& \cr
     \hline {$f_{\{2\}}$} &$0$& $\centernot{0}$ &  $0$ \cr
     \hline {$f_{\{1,2\}}$} & $\centernot{0}$ & $\centernot{0}$ &  \cr
     \hline  \end{tabular}\end{center}
Clearly, the values of $ASig(f)$, and, more precisely, the values of $ ASig_1(f)$ and $ ASig_2(f)$  are sufficient to discrimine all the possible multiple faults. This result can also be obtained by applying Criterion~\ref{Criterion3} to these two components: this criterion permits to show the equivalence between $f_1=0$ (resp. $f_2=0$) and $ ASig_{1} = 0 $  (resp. $ ASig_{2}= 0$). 

Without considering constraints on parameters and faults, the following table of signatures is obtained.
\begin{center}
\begin{tabular}{|c|c|c|c|}
\hline  
$f$  & {$ ASig_1(f)$}  & {$ ASig_2(f)$} & {$ ASig_3(f)$} \cr 
      \hline $f_{\{\}}$ &$0$&$0$&$0$\cr
     \hline {$f_{\{1\}}$} & $\centernot{0}$ &$0$& \cr
     \hline {$f_{\{2\}}$} &$0$&  &  $0$ \cr
     \hline {$f_{\{1,2\}}$} & $\centernot{0}$ &  &  \cr
     \hline  
     \end{tabular}\end{center}
This last table shows that these constraints plays an important role for studying, a priori, the values of $ASig(f)$ in function of the multiple faults. More precisely, the semialgebraic set tools insure that, for some particular values of $p_1$, $p_2$, $p_3$, $f_1$ and of $f_2$, the fault $f_{\{2\}}$ can not be detected. The same remark holds for the discrimination of the mutiple faults $f_{\{1\}}$ and $f_{\{1,2\}}$. 
\end{example}

\begin{remark}
 Algebraic criterions, using Groebner basis computations, can also be developed to obtain informations about the possible values of $ASig(f_{\mathcal{N}})$. For example, \textit{if the Groebner basis of $E_\mathcal{N} \cup  \{ ASig_k(f_{\mathcal{N}})\}$ is equal to $\{1\}$ then the $k$th component of $ASig(f_\mathcal{N})$ never vanishes.}\\
Such a criterion relies on the fact that if the sufficient condition implies that polynomials of $E_\mathcal{N} \cup  \{ ASig_k(f_{\mathcal{N}})\}$ has no common complex zeros (see~\cite{cox}) and, consequently, no real zeros. Even if constraints on parameters expressed as inequalities can not be taken into account, this criterion can be tested more rapidly in practice than Criterion~\ref{Criterion1}.
\end{remark}

\section{Application}\label{SecAppli}

The computation of the algebraic signature and the application of the three criterions had been implemented in the computer algebra system Maple 18. The table giving the expected values of the signature in function of the possible multiple faults constitutes the input of a Scilab program. This latter software is used, from the measurements of a system, to estimate numerically the algebraic signature. The comparison between the numerical values and the expected values of the signature permits to discriminate multiple  faults.\\
\\
Our method is applied on a model of two coupled water tanks (See~\cite{Nath2015,Seydou2013}) given by
\begin{equation}\label{wt1}
\left\{ \begin{array}{l}
\dot x_1(t,p)=p_1\, u(t)-p_2\, \sqrt{x_1(t,p)},\,x_1(0)=1,\\
\dot x_2(t,p)=p_3\, \sqrt{x_1(t,p)} -p_4\, \sqrt{x_2(t,p)},\,x_2(0)=0.6,\\
y(t,p)=p_5\, \sqrt{x_1(t,p)},
\end{array}\right.
\end{equation} 
where $p=(p_i)_{i=1,\ldots,5}$ is the model parameter vector, $x=(x_1,x_2)^T$ represents the state vector and corresponds to the level in each tank, and $u\not \equiv 0$ is the input vector. The water level in the tanks can vary between 0 and 10. Contrary to ~\cite{Nath2015,Seydou2013}, we suppose that there is only one output, $y$, on the first water-tank.

Let $f_1$ denote an unknown additive fault on the actuator signal, $f_2$ an additive fault on the sensor at the output of the first water tank, and $f_3\in [0;1]$ a clogging fault. The fully clogged pipe situation corresponds to $f_3=1$ and $ 0<f_3< 1$ represents a partial clogging. Afterwards, the clogged pipe situation is supposed partial.

In order to use the Rosenfled-Groebner algorithm, a change of variables is necessary. By setting $z_1(t,p)=\sqrt{x_1(t,p)}$ and $z_2(t,p)=\sqrt{x_2(t,p)}$, the model hereafter is obtained: 
\begin{equation}\label{saresoudreApplNum}
\Gamma_f
\left\{\begin{array}{l}
\dot x_1 = p_1 (u+f_1)- p_2 (1-f_3) z_1,\\
\dot x_2 = p_3 (1-f_3) z_1- p_4 z_2,\\
     z_1^2=x_1, z_2^2=x_2,\\
     y= p_5 (1-f_3) z_1+f_2
\end{array}\right.
\end{equation}
The first step of our approach can then be applied to obtained the following ARR: \\
$ 2\,y\,\dot y-p_5\,(f_3-1)^2\,(p_1\,p_5\,f_1+p_2\,f_2)-p_1\,p_5^2\,(f_3-1)^2\,u+ p_2\,p_5\,(f_3-1)^2 \,y-2\,f_2\,\dot y=0$ 
and the corresponding exhaustive summary:

$\phi(f_1,f_2,f_3) = (-p_5\,(f_3-1)^2\,(p_1\,p_5\,f_1+p_2\,f_2),\\ ~~~~~~~~~~~~~~~~~ -p_1\,p_5^2\,(f_3-1)^2 , \; p_2\,p_5\,(f_3-1)^2  ,\; -2\,f_2)\,.$\\

The second step of our method, which is the application of Algorithm \verb Algebraic_signature, provides the following signature: 

$ASig(f)= ({\phi}_1,\;{\phi}_4,\; p_1\,p_5^2+{\phi}_2,\;-p_2\,p_5+{\phi}_3,\\ ~~~~~~~~~~~~~~~~~~~~~~~~~~~~~~~ -{\phi}_3\,{\phi}_4+2\,{\phi}_1,\;-p_2\,p_5\,{\phi}_4+2\,{\phi}_1).$\\

The third step consists in computing the expected values of $ASig(f)$ in function of the multiple faults. These values are summarized in Table~\ref{tab1} (See Section~\ref{SectionCriterions} for the conventions used in this table). For this application, we start by defining the set of constraints:  $$C_{p,f} = \{0< p_1,\, \ldots, 0< p_5,\, 0 \leq f_3< 1 \} \,.$$ which corresponds to the physical signification of the model parameters and to the assumption of a non fully clogging pipe. Next, Criterions~\ref{Criterion1} and~\ref{Criterion2} are used to obtain Table~1.
 \begin{center}
\begin{small}
 \begin{tabular}{|c|c|c|c|c|c|c|} \hline 
 &&&&&&\cr
 &&&&&&\cr
 &&&&&&\cr
 &\rot{$ ASig_1(f)$}  & \rot{$ ASig_2(f)$} & \rot{$ ASig_3(f)$} & \rot{$ ASig_4(f)$} & \rot{$ ASig_5(f)$} & \rot{$ ASig_6(f)$}\cr
 \hline {$f_{\{\}}$} & 0 & 0 & 0 & 0 & 0 & 0 \cr
     \hline {$f_{\{1\}}$ } & $\centernot{0}$ & 0 & 0 & 0 & $\centernot{0}$ & $\centernot{0}$ \cr
     \hline {$f_{\{2\}}$ } & $\centernot{0}$ & $\centernot{0}$ & 0 & 0 & 0 & 0 \cr
     \hline {$f_{\{3\}}$ } & 0 & 0 & $\centernot{0}$ & $\centernot{0}$ & 0 & 0 \cr
     \hline {$f_{\{1,2\}}$} &  & $\centernot{0}$ & 0 & 0 & $\centernot{0}$ & $\centernot{0}$ \cr
     \hline {$f_{\{1,3\}}$} & $\centernot{0}$ & 0 & $\centernot{0}$ & $\centernot{0}$ & $\centernot{0}$ & $\centernot{0}$ \cr
     \hline {$f_{\{2,3\}}$} & $\centernot{0}$ & $\centernot{0}$ & $\centernot{0}$ & $\centernot{0}$ & 0 & $\centernot{0}$ \cr
     \hline {$f_{\{1,2,3\}}$} &  & $\centernot{0}$ & $\centernot{0}$ & $\centernot{0}$ & $\centernot{0}$ &  \cr
     \hline  \end{tabular}          
 
  \end{small}       
   {Table 1 \\Numerical Expected Values of the Algebraic Signatures}  
  \label{tab1}         \end{center}

Table~1 shows that the components $ ASig_2(f)$, $ ASig_4(f)$ and $ ASig_5(f)$ permit the discrimination of all the multiple faults for any input $u$. Indeed, $ ASig_2(f)$, $ ASig_4(f)$ and $ ASig_5(f)$ do not depend on the component $\phi_2$ which is the coefficient of the only term depending on $u$ in the ARR. Consequently, the model is input-strongly algebraically diagnosable.

\begin{remark}
\begin{enumerate}\label{applirem}
\item This result holds even if the values of $p_1$, $p_3$ and $p_4$ are not known. In other terms, the knowledge of the values of all the internal parameters is not needed for detecting and discriminating the possible multiple faults. 
\item The fact that the model is input-strongly algebraically diagnosable can be obtained by applying Criterion~\ref{Criterion3} to $ ASig_2(f)$, $ ASig_4(f)$ and $ ASig_5(f)$.
\end{enumerate}
\end{remark}

In the simulations, a simple controller is used to control the water level in the upper tank to follow a square reference signal. 
The parameters of the model are equal to $p_1=p_2=p_3=p_4=0.3$, $p_5=1$. 
 The simulated output are disturbed by a truncated Gaussian noise $\eta$ such that $\eta(t)\in [-0.001;0.001]$. Thus, $y(t)=\bar y(t)+\eta(t)$ where $\bar y $ is the exact output corresponding to the exact value of parameters. The observations are supposed to be done at the discrete time $(t_i)_{i=1,\ldots,M}$ on the interval $[0,50]$ with a sampling period equal to 0.5. 
In the faulty scenarios, we assume that the faults are introduced at time $t=20s$.

The derivatives are estimated in using a method based on the B-splines \cite{Ibrir}.
 In order to estimate $\phi $, the method develops in \cite{art1} is taken again. Rewriting the ARR at each discrete time $t_i$, $M$ linear relations with respect to the components of $\phi$ are obtained leading to a linear system. If we denote $y_p(t_i)$ the estimate of $\dot y(t_i)$, the system has the following form: 
 \begin{itemize}
\item in the faulty situation,
\begin{equation}\label{eqaf}
A^f X_f=b
\end{equation}
 with \\
 
 $X_f=(-p_5\,(f_3-1)^2\,(p_1\,p_5\,f_1+p_2\,f_2),\\ ~~~~~~~~~~~~~~~~~~~~  -p_1\,p_5^2\,(f_3-1)^2 , \; p_2\,p_5\,(f_3-1)^2  ,\; -2\,f_2)$, \\
$A_i^f=(1,u(t_i),y(t_i),  y_p(t_i))$ and $b_i=-2\,y(t_i)\, y_p(t_i)$.\\

Remark that $\phi(f)=(X_f(1),X_f(2),X_f(3),X_f(4))$. \\

\item in the fault-free situation,
\begin{equation}\label{eqsf}
  A X_0=b 
\end{equation}
with $X_0=(p_1\,p_5^2\,; p_2\,p_5)$,\\ $A_i=(u(t_i), y(t_i))$ and $b_i=-2\,y(t_i)\, y_p(t_i)$.
\end{itemize}


These systems will be solved with the QR factorization which does not require any initial guess.\\ 

System~(\ref{eqsf}) is used the detect the time point $t_d$ at which the multiple  fault acts. From the 10 first time points, matrix $A$ and vector $b$ are constructed. Then, at each iteration, they are completed in considering one more time and system~(\ref{eqsf}) is solved with this new matrix $A$ and this new vector $b$. The estimate of $X_0$ is compared to the nominal value obtained with the real parameter values. If their difference in norm 2 is upper than $10^{-3}$, we consider that a multiple fault acts; in other terms, the fault is algebraically detectable.\\

Once the fault detected, System~(\ref{eqaf}) serves to discriminate the multiple fault. At least four time points after the detecting time point $t_d$ are needed since $X_f$ is of dimension 4. Remark that the multiple faults can be detected and discriminated only every 0.5 second since the verification of these properties is based on the construction of systems~(\ref{eqaf}) and~(\ref{eqsf}). $ASig$ is then estimated and Table 1 is used to discriminate the multiple fault acting. The results are summarized in Table~2.\newpage

\begin{center}

\begin{tabular}{|c|c|c|}\hline
(Multi-)faults $f$ & Detection & Discrimina-\\
& times (s)& tion times (s) \\ \hline
$f_{\{1\}} =$ (\,0.5\,,\,0\,,\,0\,) &0&3\\ \hline
$f_{\{2\}} =$ (\,0\,,\,0.5\,,\,0\,)&0.5&1.5  \\ \hline
$f_{\{3\}} =$ (\,0\,,\,0\,,\,0.5\,) &0.5 &2 \\ \hline
$f_{\{1,3\}} =$ (\,0.5\,,\,0\,,\,0.1\,) &0&1.5 \\ \hline   
$f_{\{1,3\}} =$ (\,0.5\,,\,0\,,\,0.7\,)&0.5 &11$^{*}$ \\ \hline
$f_{\{1,2\}} =$ (\,0.5\,,\,0.5\,,\,0\,)&0&1.5 \\ \hline
$f_{\{2,3\}} =$ (\,0\,,\,0.5\,,\,0.1\,)&0&1.5   \\ \hline
$f_{\{2,3\}} =$ (\,0\,,\,0.5\,,\,0.7\,)&0&1.5 \\ \hline
\end{tabular}

Table 2 \\{Detection and discrimination times. }\label{tabletemps}
\end{center}
$^{*}$ $f_3 \neq 0 $ is first detected at $t=20.5s$ and the multiple  fault $f_{\{1,3\}} $ is discriminated at $t=31s$.
%
%
\section{Conclusion}\label{conclusion}
 In this paper, based on ARRs, an algebraic method for assessing (multiple)   faults diagnosability and detectability of non linear parametrized dynamical models is proposed. This method combines different algebra tools leading to efficient discriminatory relations. The application of our algorithms on the coupled water-tanks example highlights the interest of this work.

\bibliographystyle{plain}        
 \bibliography{ifacconf2}

\end{document}